\documentclass[12pt]{amsart}

\usepackage{graphicx}
\usepackage{amsmath}
\usepackage{amssymb}
\usepackage{amsfonts}
\usepackage{verbatim}
\usepackage{scalefnt}
\usepackage{multirow}
\usepackage{enumerate}
\usepackage{mathtools}
\usepackage{float}
\usepackage{enumitem}
\restylefloat{figure}
\usepackage[margin=3cm]{geometry}

\usepackage{fancyhdr}
\usepackage{hyperref}

\newtheorem{theorem}{Theorem}

\newtheorem{conjecture}[theorem]{Conjecture}

\newtheorem*{definition*}{Definition}

\newtheorem{question}[theorem]{Question}

\numberwithin{equation}{section}
\numberwithin{theorem}{section}

\title[A sumset version of a conjecture of Pilz]%
  {A sumset version of a conjecture of Pilz}

\author{J\'anos Nagy}
\email{janomo4@gmail.com}
\address{HUN-REN R\'enyi Institute of Mathematics, Budapest, Re\'altanoda u. 13-15., Budapest´1053, Hungary.;\newline \hspace*{4mm}
Department of Computer Science and Information Theory, Budapest University of Technology and Economics, M\H{u}egyetem rkp. 3., H-1111 Budapest, Hungary; \newline \hspace*{4mm}
MTA-BME Lend\"ulet Arithmetic Combinatorics Research Group, M\H{u}egyetem rkp. 3., H-1111 Budapest, Hungary.}

\author{P\'eter P\'al Pach}
\email{pachpp@renyi.hu}
\address{HUN-REN R\'enyi Institute of Mathematics, Budapest, Re\'altanoda u. 13-15., Budapest´1053, Hungary.;\newline \hspace*{4mm}Department of Computer Science and Information Theory, Budapest University of Technology and Economics, M\H{u}egyetem rkp. 3., H-1111 Budapest, Hungary; \newline \hspace*{4mm} 
MTA-BME Lend\"ulet Arithmetic Combinatorics Research Group, M\H{u}egyetem rkp. 3., H-1111 Budapest, Hungary.}

\thanks{}

\begin{document}

\begin{abstract}	
Pilz's conjecture states that for any finite set $A=\{a_1,a_2,\dots,a_k\}$ of positive integers and positive integer $n$ in the union of the sets $\{a_1,2a_1,\dots,na_1\},\dots, \{a_k,2a_k,\dots,na_k\}$ (considered as a multiset) at least $n$ values appear an odd number of times. In this short note we consider a variant of this problem. Namely, we show that in the sumset $\{a_1,2a_1,\dots,na_1\}+\dots+\{a_k,2a_k,\dots,na_k\}$ (considered as a multiset) at least $n$ values appear an odd number of times.
\end{abstract}

\date{\today}
\maketitle

\section{Introduction}

In 1992 Pilz~\cite{Pil92} formulated a conjecture about the minimal distance of a certain near-ring code. For our purposes it is convenient to formulate the conjecture in the following way:

\begin{conjecture}\label{conj-pilz}
If $n\geq 1$ and $A$ is a finite set of positive integers, then the size of the symmetric difference of the sets $A,2\cdot A,\dots,n\cdot A$ is at least $n$.     
\end{conjecture}
Here we denote by $i\cdot A=iA$ the dilation of the set $A$ by a factor $i$: 
$$i\cdot A:=\{ia:\ a\in A\}.$$
Recall that the symmetric difference $C\Delta D$ of two sets, $C,D$, is the set of elements that belong to exactly one of $C,D$, that is, $C\Delta D= (C\cup D)\setminus(C\cap D) = (C\setminus D)\cup (D\setminus C)$. Note that $\Delta$ is associative, for given sets $C_1,\dots,C_m$, their symmetric difference $C_1\Delta\dots\Delta C_m$ is simply the set of elements that belong to precisely an odd number of sets $C_i$. The particular case of Pilz’s conjecture where $A= [k]=\{1,2,\dots,k\}$ for some $k\in\mathbb{Z}^+$ was eventually established independently by Huang, Ke and Pilz \cite{HKP10} and by the second named author and C. Szabó \cite{PS11}. The general case remains open. 
There are several examples when the size of the symmetric difference is {\it exactly} $n$, for instance, when $A$ is a singleton or $A=[n]$. The currently known best lower bound for $|A\Delta (2A)\Delta\dots\Delta(nA)|$  is $\frac{n}{(\log n)^\lambda}$, where $\lambda\approx 0.2223$. \cite{PS11}
For more on Pilz's conjecture see also \cite[Section 1.4]{CCP21}

For two finite sets $A,B$ of integers let us define $A\nabla B$ to be the set of those elements that can be represented as $ab\ (a\in A,b\in B)$ in an odd number of ways. Note that for $A=\{a_1,\dots,a_k\}$ and $B=\{b_1,\dots,b_\ell\}$ we have
$$A\nabla B=(a_1B)\Delta\dots\Delta(a_kB)=(b_1A)\Delta\dots\Delta (b_\ell A).$$

By this notation Conjecture~\ref{conj-pilz} states that $A\nabla [n]$ has size at least $n$ for every finite $A\subseteq \mathbb{Z}^+$. We may switch to additive notation as follows. For finite sets $A,B\subseteq \mathbb{Z}$, or more generally, for finite subsets of the integer grid $A,B\subseteq \mathbb{Z}^r$ 
let $A\oplus B$ be the set of those elements that can be represented as $a+b\ (a\in A,b\in B)$ in an odd number of ways. Let $p_1,p_2,\dots,p_r$ denote the primes up $n$. Write each $k\leq n$ in the form $p_1^{\alpha_1}\dots p_r^{\alpha_r}$ and assign the ``exponent vector'' $v_k:=(\alpha_1,\dots,\alpha_r)$ to $k$. Let $S_k=\{v_1,v_2,\dots,v_k\}\subseteq \mathbb{Z}^r=\mathbb{Z}^{\pi(n)}$. For instance, in case of $n=4$ we get the $L$-shape $S_4=\{(0,0),(1,0),(0,1),(2,0)\}\subseteq \mathbb{Z}^2$. 

By this notation Conjecture~\ref{conj-pilz} states that $|S_n\oplus A|\geq n$ for any finite set $A\subseteq\mathbb{Z}^{\pi(n)}$. Alternatively, the conjecture states that the symmetric difference of finitely many translates of $S_n$ always has size at least $n$. It is a nice exercise to show that $|S\oplus A|\geq |S|$ holds if $S=\{0,1\}^r$ is a 2-cube, Pilz's conjecture states that $S=S_n$ also satisfies this inequality. 

However, in general, the inequality $|S\oplus A|\geq |S|$ may not hold, already in dimension 1, the set $S\oplus A$ can be much smaller than $S$. For instance, for $S=[n]$ and $A=\{0,1\}$ we get that the set $S\oplus A=\{1,n+1\}$ has only two elements. For the inequality to hold we shall require further conditions on the sets $S$ and $A$. In this note we consider the following (1-dimensional) sumset variant of Pilz's conjecture:

\begin{question}\label{ques-1}
Is is true that $|S_1\oplus S_2\oplus \dots \oplus S_k|\geq n$, if each $S_i$ is of the form $S_i=\{a_i,2a_i,\dots,na_i\}$ for some $a_i\in\mathbb{Z}^+$?
\end{question}
\noindent
We answer this question in the affirmative:
\begin{theorem}\label{thm-1}
Let $n,k$ be  positive integers. If $a_1,a_2,\dots,a_k\in\mathbb{Z}^+$, then 
$$|{\oplus}_{i=1}^k \{a_i,2a_i,\dots,na_i\}|\geq n.$$
\end{theorem}
\noindent
In fact we prove a slightly stronger statement:
\begin{theorem}\label{thm-2}
Let $n,k$ be  positive integers and $V\subseteq \mathbb{Z}^+$ a finite set of odd size. \\
If $a_1,a_2,\dots,a_k\in\mathbb{Z}^+$, then 
$$|V{\oplus} {\oplus}_{i=1}^k \{a_i,2a_i,\dots,na_i\}|\geq n.$$
\end{theorem}

Finally, we shall mention an open problem from geometry of similar nature: Is it true that the area of the symmetric difference of an odd number of unit discs is always at least $\pi$? This was first asked by Pak \cite{Pak}, the problem is still open, for more on this problem, see also \cite{Pin23}. 
However, there is an important difference between this problem from combinatorial geometry and Pilz's conjecture (and Question~\ref{ques-1}): in case of the latter problems it is {\it not} assumed that we take an {\it odd} number of translates of the corresponding set (the size of $A$ in Pilz's conjecture and the number $n$ in Question~\ref{ques-1} may be even).

\section{Proof of Theorem~\ref{thm-1}~and~\ref{thm-2}}

Let us assign a polynomial $p_S(x)\in\mathbb{F}_2[x]$ to each finite subset $S$ of nonnegative integers: $p_S(x)=\sum\limits_{s\in S} x^s$. Observe that $p_{S_1\Delta S_2}(x)=p_{S_1}(x)+p_{S_2}(x)$ and $p_{{\oplus}_{i=1}^k S_i}(x)=\prod_{i=1}^k p_{S_i}(x) $.  

Therefore, proving Theorem~\ref{thm-1} is equivalent to showing that the number of nonzero coefficients in $\prod\limits_{i=1}^k \left(x^{a_i}+x^{2a_i}+\dots+x^{na_i} \right)$ is at least $n$. After expanding out $\prod_{i=1}^k x^{a_i}$ we get the polynomial
$$p(x):=\prod\limits_{i=1}^k \left(1+x^{a_i}+x^{2a_i}+\dots+x^{(n-1)a_i} \right),$$
our aim is to show that the number of nonzero coefficients of $p$ is at least $n$. 
Without loss of generality, we may assume that $\gcd(a_1,\dots,a_k)=1$, since otherwise we may consider $p(x)$ as a polynomial of  $x^{\gcd(a_1,\dots,a_k)}$ which has the same number of nonzero coefficients as $p$.

Let $n=2^\alpha t$, where $\alpha\geq 0$ and $t$ is odd. 

Let us write $p(x)$ as $p(x)=q(x)r(x)$, where 
$$q(x)=\prod\limits_{i=1}^k \left(1+x^{a_i}+x^{2a_i}+\dots+x^{(t-1)a_i} \right),$$
$$r(x)=\prod\limits_{i=1}^k \left(1+x^{ta_i}+x^{2ta_i}+\dots+x^{(2^\alpha-1)ta_i} \right).$$

First, we turn our attention at $q(x)$. Let us write $q(x)$ as \begin{equation}\label{qdef}
q(x)=q_0(x^t)+xq_1(x^t)+\dots+x^{t-1}q_{t-1}(x^t),
\end{equation}
that is, we partition the monomials in $q$ into $t$ groups according to the mod $t$ residue of the exponent of $x$.  
Then $$p(x)=q(x)r(x)=q_0(x^t)r(x)+xq_1(x^t)r(x)+\dots+x^{t-1}q_{t-1}(x^t)r(x),$$
where the nonzero coefficients of these $t$ polynomials are pairwise different, since $r(x)$ is also a polynomial of $x^t$. Hence, it suffices to prove that each $q_i(x^t)r(x)$ has at least $2^{\alpha}$ nonzero coefficients. 

We show that $q_i(1)=1$ for each $i$, that is, the number of nonzero coefficients of $q_i$ is odd.

If we expand out $q(x)$, then the number of terms -- without cancellations -- is $t^k$, which is odd. We show that their exponents are uniformly distributed modulo $t$, implying that each residue is obtained $t^{k-1}$ times, thus $q_i(1)=1$ indeed holds. For a residue $b$ (modulo $t$) let $F(b)$ denote the number of terms (before cancellations) where the exponent has residue $b$ mod $t$. Since $(1+x^{a_i}+\dots +x^{(t-1)a_i})$ is among the factors, $F(b)=F(b+a_i)$ for every $b$. This holds for every $i$ and the greatest common divisor of the numbers $a_1,\dots,a_k$ is 1, so $F$ is constant. Thus $q_i(1)=1$, as we claimed.

Now, we show that $q_i(x^t)r(x)$ has at least $2^\alpha$ nonzero coefficients.

Setting $y=x^t$ we have
$$r(x)=\prod\limits_{i=1}^k \left(1+y^{a_i}+y^{2a_i}+\dots+y^{(2^\alpha-1)a_i} \right)=\prod\limits_{i=1}^k \left(1+y^{a_i}\right)^{2^\alpha-1}.$$
Let $a_i=2^{\alpha_i}t_i$, where $\alpha_i\geq 0$ and $t_i$ is odd. By using the identity
$$1+y^{a_i}=(1+y^{t_i})^{2^{\alpha_i}}=(1+y)^{2^{\alpha_i}}(1+y+\dots+y^{t_i-1})^{2^{\alpha_i}},$$
we get that
$$r(x)=(1+y)^{(2^\alpha-1)\sum_{i=1}^k 2^{\alpha_i}} \prod\limits_{i=1}^k (1+y+\dots+y^{t_i-1})^{(2^{\alpha}-1)2^{\alpha_i}}.$$
Let us express the exponent of $1+y$ as a sum of distinct 2-powers: 
$$(2^\alpha-1)\sum_{i=1}^k 2^{\alpha_i}=\sum\limits_{j\in J} 2^{\beta_j}.$$
Note that $|J|\geq \alpha$. (Indeed, the possible residues that a 2-power can have modulo $2^{\alpha}-1$ are $1,2,2^2,\dots,2^{\alpha-1}$. Assume we get the 0 residue with a sum containing a minimum number of terms. Then all the residues are distinct in the sum, since otherwise two copies of a 2-power $2^j$ may be replaced by one copy of $2^{j+1}$, which would contradict minimality. However, if all terms are distinct, then we have to add all of them to get 0 mod $2^{\alpha}-1$, since their total sum is exactly $2^{\alpha}-1$.) 

Consider the set $S:=\{ \sum\limits_{j\in J} \varepsilon_j 2^{\beta_j}:\ \varepsilon_j\in \{0,1\} \}$ and observe that $$(1+y)^{(2^\alpha-1)\sum_{i=1}^k 2^{\alpha_i}}=(1+y)^{\sum_{j\in J}2^{\beta_j}}=\sum_{s\in S}y^s=p_S(y).$$
If we write $$q_i(y)\prod\limits_{i=1}^k (1+y+\dots+y^{t_i-1})^{(2^\alpha-1)2^{\alpha_i}}=\sum_{u\in U}y^u,$$
then $r(x)q_i(y)=\sum_{s\in S}\sum_{u\in U} y^{s+u}$. Note that $|U|$ is odd, since $q_i(1)=1$ and each $t_i$ is odd. For estimating the number of nonzero coefficients of $r$ we shall compute the size of the symmetric difference of the sets $S+u=\{s+u:s\in S\}$ ($u\in U$). Our aim is to show that the size of this symmetric difference is at least $|S|$. 
To see this, we first prove that $S$ tiles $\mathbb{Z}_{\geq 0}$, that is, there is some $R$ such that $\mathbb{Z}_{\geq 0}$ is the direct sum of $S$ and $R$, meaning that every nonnegative integer can be uniquely represented as $s+r$ with $s\in S,r\in R$. This is immediate, since we can choose $R$ to be the set of those nonnegative integers whose base-2 representation does not contain any of $2^{\beta_j}$ ($j\in J$). 
Now, we show that there is an $|S|$-colouring of $\mathbb{Z}_{\geq 0}$ such that each translate $S+u$ contains exactly one element from each colour class. 
If $m\geq 0$ is an integer, then $m$ can be uniquely written as $m=s+r$ with $s\in S, r\in R$. Let us define the colour of $m$ to be $s$. 
Let us consider a translate $S+v$ and assume that $s'+v$ and $s''+v$ have the same colour. Then $s'+v=s+r_1$ and $s''+v=s+r_2$ for some $s\in S, r_1,r_2\in R$. However, these equations imply that $s'+r_2=s''+r_1$, but $S+R$ is a direct sum, so $s'=s''$ and $r_2=r_1$. Therefore, the colouring satisfies our requirements.

Since $|U|$ is odd, the symmetric difference of the sets $S+u$ contains an odd number of elements from each of the $|S|$ colour classes, thus its size is indeed at least $|S|=2^{|J|}\geq 2^\alpha$. 

Hence, $q_i(x^t)r(x)$ has at least $2^\alpha$ terms. This concludes the proof of Theorem~\ref{thm-1}.

Theorem~\ref{thm-2} can be proved similarly: The only difference is that in \eqref{qdef} we shall write $q(x)\sum\limits_{v\in V} x^v$ in place of $q(x)$ on the left hand-side of the equation. When we expand out $q$, the exponents are uniformly distributed modulo $t$ (before the cancellations), so the same holds for the exponents of the terms arising in $q(x)\sum\limits_{v\in V} x^v$. Since $|V|$ is odd, the rest of the argument is also fine in this setting.

\section{Acknowledgements}
Both authors were supported by the Lend\"ulet program of the Hungarian Academy of Sciences (MTA). PPP was also supported by the National Research, Development and Innovation Office NKFIH (Grant Nr. K146387).

\medskip

\end{document}